\documentclass{article}

\usepackage{amsmath, amssymb, latexsym}

\usepackage{amsthm}

\newtheorem{theorem}{Theorem}[section]
\newtheorem{proposition}[theorem]{Proposition}
\newtheorem{lemma}[theorem]{Lemma}

\theoremstyle{definition}
\newtheorem{definition}[theorem]{Definition}
\theoremstyle{remark}
\newtheorem{remark}[theorem]{Remark}
\theoremstyle{remark}

\theoremstyle{remark}

\begin{document}

\title{Full Groups and Orbit Equivalence in Cantor Dynamics}

\author{{\bf K. Medynets}
 \\
 Ohio State University,  Columbus, USA
 \\
 medynets@math.ohio-state.edu
 }

 \date{}

 \maketitle

\begin{abstract} In this note we consider dynamical systems $(X,G)$ on a Cantor set $X$ satisfying some mild technical conditions.
The considered class includes, in particular, minimal and transitive
aperiodic systems. We prove that two such systems $(X_1,G_1)$ and
$(X_2,G_2)$ are orbit equivalent if and only if their full groups
are isomorphic as abstract groups. This result is a topological
version of the well-known Dye's theorem established originally for
ergodic measure-preserving actions.
\end{abstract}

\section{Introduction} Denote by $Homeo(X)$ the group of all homeomorphisms of a Cantor set
$X$.  Then the pair $(X,G)$, where $G$ is a subgroup of $Homeo(X)$,
is called a {\it Cantor dynamical system}. We would like to
emphasize from the very beginning that, in contrast to the ergodic
theory, the group $G$ does not have to be countable.

For a point $x\in X$, denote by $Orb_G(x) = \{g(x) | g\in G\}$ its
{\it $G$-orbit}. Set also $$[G] = \{\gamma \in Homeo(X) : \gamma(x)
\in Orb_G(x)\mbox{ for all }x\in X\}.$$ Then $[G]$ is a subgroup of
$Homeo(X)$, which is called the {\it full group of $G$}.

Two dynamical systems $(X,G)$ and $(Y,H)$ on the Cantor sets $X$ and
$Y$  are called {\it orbit equivalent} if there is a homeomorphism
$\Lambda : X\rightarrow Y$ such that $\Lambda(Orb_G(x)) =
Orb_H(\Lambda(x))$ for all $x\in X$.

The notion of orbit equivalence and full groups appeared first in
the context of ergodic theory in the seminal papers \cite{Dye1} and
\cite{Dye2}. In these works Dye established a number of remarkable
results that classified measure-preserving dynamical systems up to
orbit equivalence. In particular, he showed the following result.
Let two countable groups $G$ and $H$ act on a standard measure space
$(Y,\mu)$ by measure-preserving automorphisms. Then the full groups
$[G]$ and $[H]$ (defined by measure-preserving transformations) are
isomorphic as abstract groups if and only if the systems $(Y,\mu,G)$
and $(Y,\mu,H)$ are orbit equivalent. Furthermore, the algebraic
isomorphism between full groups is always spatially generated by a
map that implements the orbit equivalence.

This result implies  that the full group ``remembers'' all the
dynamical information which does not depend on the order of points
within orbits. To the best of our knowledge there is still no a
complete algebraic description of full groups. However,  some
partial results clarifying how different dynamical properties affect
the algebraic structure of the full group have been earlier
established, see, for example, \cite{eigen:1981},
\cite{KitTsan:2010}, and \cite{Mercer:1993}.

We should point out that Dye's theorem is a universal result as it
holds in completely different dynamical setups. For example, there
is a Borel version of Dye's theorem \cite{MilRosen:2007} established
for full groups of Borel equivalence relations. The thesis of Miller
\cite{Miller:Thesis} contains  algebraic characterizations (in terms
of full groups) of certain properties of underlying Borel dynamical
systems.

In the context of Cantor dynamics, the topological version of Dye's
theorem was earlier obtained for  minimal actions of locally finite
groups and the group $\mathbb Z$, see \cite{gps:1999}. In view of
the work \cite{gpsm:2010}, this result is expandable to minimal
$\mathbb Z^n$-actions. Some algebraic properties of the full group
$[\mathbb Z]$ are present in the papers \cite{Matui:2006} and
\cite{BezuglyiMedynets:2008}. We should also mention the work
\cite{Matui:2002}, where a version of Dye's theorem is established
for minimal actions of the group $\mathbb Z$ on locally compact
zero-dimensional spaces.

The main result of the present paper (Theorem
\ref{TheoremMainTheorem}) is the extension of Dye's theorem on
almost arbitrary Cantor systems. Namely, we prove that two Cantor
systems $(X,G)$ and $(Y,H)$, which meet some mild technical
conditions, are orbit equivalent if and only if their full groups
are isomorphic. This result shows that  hyperfinite and
non-hyperfinite actions can be already distinguished
 at the level of full groups (cf. \cite{KitTsan:2010} for the ergodic
 case).

 After the paper was submitted, we became aware of the work
 \cite{Rubin:book1989} devoted to the reconstruction of Boolean algebras from their transformation
 groups. In particular, Theorem 4.5(c) in there implies the main result of the present paper
 (after the corresponding interpretation of the result).  It can be derived  from \cite[Theorem
4.5(c)]{Rubin:book1989} that if two Cantor systems have orbits at
least of length three and the set of points with orbits of length
six is nowhere dense, then any isomorphism between full groups is
spatially generated. We would like to mention, though, that our
proof is significantly different and shorter from that of
\cite[Theorem 4.5(c)]{Rubin:book1989} and requires less
prerequisites.

\section{Spatial Realization} In this section we establish the main
result of the paper.  In our proof we will use a result of Fremlin
\cite[Theorem 384D]{fremlin} that states that algebraic isomorphisms
between groups of automorphisms of complete Boolean algebras are
always generated by an automorphism of the underlying algebras. We
will apply this result to the case of full groups and show that the
automorphism gives rise to a homeomorphism of the Cantor sets,
which, in its turn, implements an orbit equivalence. The method of
\cite[Theorem 384D]{fremlin} has been already used in Cantor
dynamics (see \cite{BezuglyiMedynets:2008}) to show that the
commutator of the topological full group of  minimal $\mathbb
Z$-action is a complete invariant for flip conjugacy. It should be
noted that the Boolean algebra automorphism obtained in
\cite{BezuglyiMedynets:2008}  automatically turned out to be a
homeomorphism (due to the definition of the topological full group).
In the general case, we have to find an algebraic criterion for a
set to be clopen. To achieve our program, we will need some notions
of the theory of Boolean algebras.

Let $X$ be a Cantor set. Recall that an open set $A$ is called {\it
regular open} if $A=int(\overline{A})$. Denote the family of all
regular open sets by $RO(X)$. Notice that the family of clopen sets
(denoted by $CO(X)$) is contained in $RO(X)$.

Let $\mathcal A$  be a Boolean algebra and $\mathcal H\subset
\mathcal A$. Define $sup(\mathcal H)$ to be the smallest element of
$\mathcal A$ that contains all elements of $\mathcal H$. If
$sup(\mathcal H)$ exists for any family $\mathcal H\subseteq
\mathcal A$, then the Boolean algebra $\mathcal A$ is called {\it
complete}.

\begin{proposition}
$RO(X)$ is a complete Boolean algebra with Boolean operations given
by
$$A + B= int (\overline{A\cup B}),\;\; A\cdot B= A\cap B,\;\;
A - B=A\setminus \overline{B}$$ and with suprema given by
$sup(\mathcal H)=int(\overline{\bigcup \mathcal H })$.
\end{proposition}
{\it Proof.} See  Theorem 314P of \cite{fremlin}. \hfill$\square$

\begin{remark} Notice that  finite set-theoretical unions, intersections, and
complements of clopen sets coincide with the corresponding Boolean
operations.
\end{remark}

 For a homeomorphism $\gamma$, define its {\it support} as the set
 $$spr(\gamma) = int\overline{\{x\in X | \gamma(x)\neq x\}}.$$
Note that $spr(\gamma)$ is a regular open set. Since every
homeomorphism $\gamma$ defines a Boolean algebra $RO(X)$
isomorphism, we can also define the support of $\gamma$ as the least
regular open set $S(\gamma)$ such that $\gamma(V) = V$ for every
regular open set $V\subset X-S(\gamma)$ (see  \cite[Def.
381B]{fremlin}). It is not hard to check that $S(\gamma) =
spr(\gamma)$.

We will also consider a point-wise support of the homeomorphism
$\gamma$ defined by $supp(\gamma) = \{x\in X | \gamma(x)\neq x\}$.
We note that both sets $supp(\gamma)$ and $spr(\gamma)$ are
$\gamma$-invariant and open. Furthermore,  $supp(\gamma)\subset
spr(\gamma)$ and for any point $x\in supp(\gamma)$ there is a clopen
neighborhood $W$ such that $\gamma(W)\cap W=\emptyset$.

\begin{definition} Following \cite[Def. 382O]{fremlin}, we say that a  group
$\Gamma\subset Homeo(X)$  {\it has many involutions} if for any
regular open set $A$ there is an involution $\pi\in \Gamma$ with
$spr(\pi)\subset A$.
\end{definition}

\begin{proposition} Let $(X,G)$ be a Cantor dynamical system.

(1) If for every clopen set $A$ there is a point $x\in A$ whose
$G$-orbit intersects $A$ at least twice, then then the full group
$[G]$ has many involutions.

(2) If every orbit of $G$ has the length at least three, then the
supports of involutions (with clopen supports) from $[G]$ generate
the Boolean algebra $CO(X)$ (with the standard set-theoretical
operations).
\end{proposition}
{\it Proof.} (1) For a clopen set $A$, find a point $x\in A$ and
$g\in G$ such that  $g(x)\in A$ and $x\neq g(x)$. Choose a clopen
set $V\ni x$ with $V,g(V)\subset A$ and $V\cap g(V)=\emptyset$.
Define $\pi|V=g|V$, $\pi|g(V) = g^{-1}|g(V)$, and $\pi = id$
elsewhere. Then $\pi$ is an involution supported by $A$.

(2) Fix a clopen set $A$ and a point $x\in A$. If there is an
element $g_x\in G$ with $g_x(x)\in A$ and $g_x(x)\neq x$, then
choose a clopen set $U_x$ such that $g_x(U_x)\cap U_x = \emptyset$
and $U_x,g_x(U_x)\subset A$. If for all $g\in G\setminus \{id\}$,
$g(x)\notin A$, then we choose two elements $g_x^{(1)},g_x^{(2)}\in
G$ with $\{x,g_x^{(1)}(x), g_x^{(2)}(x)\}$ being distinct points.
Choose a set $U_x\subset A$ so that
$\{U_x,g_x^{(1)}(U_x),g_x^{(2)}(U_x)\}$ are mutually disjoint and
$g_x^{(1)}(U_x)$, $g_x^{(2)}(U_x)$ are subsets of $X\setminus A$.
Take a finite subcover $\{U_{x_1},\ldots,U_{x_n}\}$ of $A$. If
$g_{x_i}(U_{x_i})\subset A$, then we may construct an involution
$\pi_i$ as in (1) with $x\in spr(\pi)$ and $spr(\pi)$ being a clopen
subset of $A$. If the trajectory of $x$ intersects $A$ only once,
then by using the elements $g_{x_i}^{(1)}$ and $g_{x_i}^{(2)}$, we
may construct two involutions $\pi_i^{(1)}$ and $\pi_i^{(2)}$ as in
(1) so that $spr(\pi_i^{(1)})$ and $spr(\pi_i^{(2)})$ are clopen
sets and $spr(\pi_i^{(1)})\cap spr(\pi_i^{(2)}) = U_{x_i}$. This
implies that the set $A$ is represented as a union and intersection
of a finite number of clopen supports of involutions.
\hfill$\square$

\medskip As a corollary, we get that every  transitive system $(X,G)$ with infinite orbits has many involutions.
 Furthermore, the clopen supports of the involutions from $[G]$ generate $CO(X)$.

\begin{theorem}\label{TheoremMainTheorem} Let $(X_1,G_1)$ and $(X_2,G_2)$ be Cantor dynamical
systems such that   each $G_i$-orbit contains at least three points
and the full group $[G_i]$ has many involutions for $i=1,2$. Then
$(X_1,G_1)$ and $(X_2,G_2)$ are orbit equivalent if and only if
 $[G_1]$ and $[G_2]$ are isomorphic as abstract groups.

Furthermore, for  every isomorphism $\alpha:[G_1]\rightarrow [G_2]$
there is a  homeomorphism $\Lambda : X_1\rightarrow X_2$ such that
$\alpha(g) = \Lambda
 g\Lambda^{-1}$ for all $g\in [G_1]$.
\end{theorem}
{\it Proof.} It is obvious that the orbit equivalence implies the
isomorphism of full groups. Conversely, let $\alpha:[G_1]\rightarrow
[G_2]$ be a group isomorphism. Since both groups $[G_1]$ and $[G_2]$
have many involutions, there is an automorphism of Boolean algebras
$\Lambda: RO(X_1)\rightarrow RO(X_2)$ such that $\alpha(g)(V) =
\Lambda g\Lambda^{-1}(V)$ for any $g\in [G_1]$ and any $V\in
RO(X_2)$ (see Theorem 384D in \cite{fremlin}). Our goal is to show
that $\Lambda$ gives rise to a homeomorphism of the Cantor sets by
establishing that $\Lambda(CO(X_1)) = CO(X_2)$.

Let $G$ stand for either of the groups $G_1$ and $G_2$. We will
establish some general properties of the full group $[G]$. For any
regular open set $V$, set $\Gamma_V = \{\gamma\in [G] :
spr(\gamma)\subset V\}$. It will be clear from the context which
group $G_1$ or $G_2$ is meant.  The following lemma immediately
follows from the proof of Theorem 384D of \cite{fremlin} (see items
(g) and (i) therein). We also refer the reader to the proof of
Theorem 5.8 in \cite{BezuglyiMedynets:2008} to see how this result
can be obtained from the scratch in the case of topological full
groups of minimal $\mathbb Z$-systems.
\begin{lemma}\label{LemmaDefAutomorphism}
(1) Let $V\in RO(X_1)$. Then  $\Lambda(V) = sup\{spr(\alpha(\pi))  :
\pi\in [G_1]\mbox{ is an involution with }spr(\pi)\subset V\}$.

(2)  If $\pi\in [G_1]$ is an involution, then $spr(\alpha(\pi)) =
\Lambda(spr(\pi))$ and $$\alpha(\Gamma_{spr(\pi)}) =
\Gamma_{spr(\alpha(\pi))}.$$
\end{lemma}

\begin{definition} For an involution $\pi\in [G]$, set $W_\pi =
\Gamma_{spr(\pi)}$. Then, in view of (2) in Lemma
\ref{LemmaDefAutomorphism},  $\alpha(W_\pi) = W_{\alpha(\pi)}$. We
note that the proof of \cite[Theorem 384D]{fremlin} contains a
precise algebraic description of the subgroups $W_\pi$. See also
Corollary 2.10 in \cite{Rubin:1996}.

\end{definition}

\medbreak Our goal now is to give an algebraic criterion for
$spr(\pi)$ to be a clopen set. We will need the following two
lemmas. For any group $\Gamma$ in $[G]$ denote by $\Gamma^\perp$ the
commutator of $\Gamma$ in $[G]$, i.e. $\Gamma^\perp = \{\rho\in [G]
: \rho \gamma = \gamma \rho\mbox{ for any }\gamma\in\Gamma\}$. For
every set $V\in RO(X)$, set also $V^\perp = X -  V$.

%
%
\begin{lemma}\label{LemmaCommutatorDescription} If  $V$ is a regular open set,
then $$\Gamma_V^\perp = \{\gamma\in [G] : spr(\gamma)\subset X - V\}
= \Gamma_{V^\perp}.$$
\end{lemma}
{\it Proof.} Suppose that $\gamma\in[G]$ such that $spr(\gamma)
\subset X - V$. Fix any element $\rho\in \Gamma_V$. Then
$spr(\rho)\subset V$. Hence $spr(\rho)\cap spr(\gamma) = \emptyset$.
This implies that $\rho$ and $\gamma$ commute.

Conversely, if $\gamma \in \Gamma_V^\perp$ and
$spr(\gamma)\not\subset X-V = X\setminus \overline V$, then
$spr(\gamma)\cap \overline V \neq \emptyset$. Since $spr(\gamma)$ is
an open set, we get that $spr(\gamma)\cap V \neq \emptyset$.   Since
the pointwise support $supp(\gamma)$ of $\gamma$  is an open dense
subset of $spr(\gamma)$, we get that $supp(\gamma)\cap
V\neq\emptyset$.

Take  any clopen subset $W$  of $supp(\gamma)\cap V$ such that
$\gamma(W)\cap W = \emptyset$. Using the fact that the group $[G]$
has many involutions, find an involution $\rho$ supported by $W$.
Clearly, $\rho\in \Gamma_V$. Take any clopen set $O\subset W$ with
$\rho(O)\cap O = \emptyset$. Then $\rho\gamma(O)\cap
\gamma\rho(O)=\emptyset$, which is a contradiction. \hfill$\square$

\medskip For a set $F\subset [G]$, denote by $<F>$ the subgroup
generated by the elements of $F$.

\begin{lemma} A regular open set $V$ is clopen if and only if for any involution $\pi\notin <\Gamma_V,\Gamma_{V^\perp}>$ there is an element
$\rho\in <\Gamma_V,\Gamma_{V^\perp}>$ such that $h = \pi\rho$ is an
involution and

\begin{itemize}
\item[(i)] if  $g\in \Gamma_V\cap W_h$, then $h^{-1}gh\in \Gamma_{V^\perp}$;

\item[(ii)] if  $g\in \Gamma_{V^\perp}\cap W_h$, then $h^{-1}gh\in
\Gamma_{V}$.

\end{itemize}

\end{lemma}
{\it Proof.} (1) First of all assume that $V$ is clopen. Set  $R =
<\Gamma_V,\Gamma_{V^\perp}>$. Let $\pi\notin R$. Set $$A =
V\cap\pi(V)\cap spr(\pi)\mbox{ and }B = V\cap \pi(X\setminus V)\cap
spr(\pi).$$ Since the set $V$ is clopen,  both of the sets $A$ and
$B$ are  open. Furthermore, $B\neq\emptyset$ as $\pi\notin R$ and
$\pi(A) = A$. Observe that $\overline A\cap\overline B = \emptyset$.
Indeed, if otherwise, take $x\in \overline A\cap\overline B$, then
there are two sequences $\{a_n\}\subset A$ and $\{b_n\}\subset B$
with $a_n\to x$ and $b_n\to x$. However, by the definition of $A$
and $B$, we get that $\{\pi(a_n)\}\subset V$ and
$\{\pi(b_n)\}\subset X\setminus V$. As the set $V$ is clopen, we get
that $\pi$ is not continuous at $x$, which is a contradiction.

Thus, we can choose a clopen set $O\subset V$ such that $A\subset O$
and $O\cap \overline B = \emptyset$. Note that $\pi(O) = O$.  Set
$\rho_1|O = \pi^{-1}|O$ and $\rho_1 = id$ elsewhere. Then
$\pi\rho_1|(V\setminus B) = id$. Clearly, $\pi\rho_1$ is an
involution and $\rho_1\in \Gamma_V$.

Repeating the same arguments with the set $X\setminus V$ and the
involution $\pi\rho_1$, we find an element $\rho_2\in
\Gamma_{V^\perp}$ such that $\pi\rho_1\rho_2|(V^\perp\setminus
\pi(B)) = id$. Set $\rho = \rho_1\rho_2\in R$ and $h=\pi\rho$.
Observe that $spr(h) = B\cup \pi(B)$ and $h|B = \pi | B$.

Now if $g\in W_h\cap \Gamma_V$, then $spr(g)\subset B$. Hence
$spr(hgh^{-1})\subset \pi(B)$ and $hgh^{-1}\in \Gamma_{V^\perp}$.
The condition (ii) can be established in a similar way.

(2) Conversely, assume that $V$ is a non-closed regular open set.
Set
$$B = X\setminus(V\cup(X-V)).$$ Then $B$ is a non-empty closed set.
Note also that $B = \overline V\cap \overline {X-V}$.

(2-i) Fix a point $x\in B$. If there is an element $g\in [G]$ with
$g(x)\in V$ or $g(x)\in X-V$. Take a small clopen neighborhood $U$
of $x$ with $g(U)\cap U = \emptyset$. Define an involution
$\pi\notin R$ by setting $\pi|U = g|U$, $\pi|g(U) = g^{-1}|g(U)$,
and $\pi = id$ elsewhere.

Without loss of generality, we will assume that $g(x) \in V$. Fix
any element $\rho\in R$. Since $\rho(x) = x$, we get that
$\pi\rho(x)\neq x$. Set $h = \pi\rho$. Hence, $h(O)\cap O =
\emptyset$ for some clopen neighborhood of $x$. Note that $O\cap V
\neq \emptyset$. Take any involution $q\in \Gamma_V$ with
$spr(q)\subset O$. It follows that $q\in W_{h}$ and
$spr(hqh^{-1})\subset h(O)\subset V$. Hence $hqh^{-1}\in \Gamma_V$,
which contradictions the condition (i) of the lemma. The case when
$g(x)\in X -  V$ is proved similarly.

(2-ii). Now consider the situation when $g(x)\in B$ for any point
$x\in B$ and any element $g\in [G]$. This means that $g(B) = B$ for
all $g\in G$. Since, every $G$-orbit has the length at least three,
choose three different points $\{x_1,x_2,x_3\}\subset B$ from the
same $G$-orbit. By the standard arguments, find two involutions
$\pi_1$ and $\pi_2$ from $[G]$ with clopen supports such that
$\pi_1(x_1) = x_2$ and $\pi_2(x_2) = x_3$.

 We claim that there is a point $x\in\{x_1,x_2,x_3\}$ and a
homeomorphism $\pi\in\{\pi_1,\pi_2, \pi_1\pi_2\}$ such that for
every clopen neighborhood $U$ of $x$ either $\pi(U\cap V)\cap V\neq
\emptyset$ or $\pi(U\cap V^\perp)\cap V^\perp\neq \emptyset$.

Assume the converse. Since $\pi_1(B) = B$, there is a clopen
neighborhood $U_1$ of $x_1$ with $\pi_1(U_1\cap V)\subset V^\perp$
and $\pi_1(U_1\cap V^\perp)\subset V$. In the same way, there is a
clopen neighborhood $U_2$ of $x_2$ with $\pi_2(U_2\cap V)\subset
V^\perp$ and $\pi_2(U_2\cap V^\perp)\subset V$. It follows that
$\pi_2\pi_1(U\cap V)\subset V$ for some clopen subset $U\ni x$ of
$U_1$, which is a contradiction.

Without loss of generality, we will assume that $\pi(U\cap V)\cap
V\neq \emptyset$ for every clopen neighborhood of $x$.  This means
that every open neighborhood of $x$ contains an element $g\in
\Gamma_V$ with $\pi g\pi^{-1}\in \Gamma_V$. Repeating arguments from
(2-i), we get that for any $\rho\in R$ there is $g\in \Gamma_V$ with
$(\pi\rho) g (\pi\rho)^{-1}\in \Gamma_V$.

This proves the necessity of the condition (i). The necessity of the
condition (ii) can be established in a similar way. \hfill$\square$

\medskip {\it Continuation of the proof}. It follows from the lemmas
above that if a clopen set $V$ is the support of an involution
$\pi$, then $\Lambda(V) = spr(\alpha(\pi))$ is a clopen set as well.

Since $CO(X_i)$ is generated by the  clopen supports of involutions
 and $\Lambda$ is a Boolean algebra isomorphism,
we conclude that $\Lambda(CO(X_1)) = CO(X_2)$. This implies that
$\Lambda$ defines a homeomorphism of $X_1$ and $X_2$.

If $ y = g(x)$ for some $x\in X_1$ and $g\in [G_1]$, then
$\Lambda(y) = \Lambda g(x) = \alpha(g) \Lambda(x)$. Thus,
$\Lambda(x)$ and $\Lambda(y)$ belong to the same $G_2$-orbit. Hence,
$\Lambda$ implements the orbit equivalence between $(X_1,G_1)$ and
$(X_2,G_2)$. \hfill$\square$

\begin{definition} Let $(X,G)$ be a Cantor dynamical system. Then
the {\it topological full group} of $G$ (in symbols $[[G]]$) is
formed by all elements $\gamma\in [G]$ for which there is a clopen
partition $\{U_1,\ldots, U_n\}$ of $X$ and elements $g_1,\ldots,
g_n\in G$ such that $\gamma|U_i = g_i|U_i$ for every $i$.
\end{definition}

\begin{remark} Note that if $[G]$ has many involutions, then so does the topological full group $[[G]]$ (see \cite{gps:1999} for the
definition). So we can follow the lines of the proof of Theorem
\ref{TheoremMainTheorem} to  show that any group isomorphism between
topological full groups is  spatially generated.
\end{remark}

\medskip {\bf Acknowledgement.} I would like to thank Sergey
Bezuglyi for introducing this subject to me and for numerous helpful
discussions. I am also thankful to Matatyahu Rubin for the
discussions of reconstruction theorems.

\end{document}